\documentclass[12pt,a4paper]{amsart}
\usepackage{mathrsfs}

\newtheorem{theo+}              {Theorem}           [section]
\newtheorem{prop+}  [theo+]     {Proposition}
\newtheorem{coro+}  [theo+]     {Corollary}
\newtheorem{lemm+}  [theo+]     {Lemma}
\newtheorem{exam+}  [theo+]     {Example}
\newtheorem{rema+}  [theo+]     {Remark}
\newtheorem{defi+}  [theo+]     {Definition}

\newenvironment{theorem}{\begin{theo+}}{\end{theo+}}
\newenvironment{proposition}{\begin{prop+}}{\end{prop+}}
\newenvironment{corollary}{\begin{coro+}}{\end{coro+}}
\newenvironment{lemma}{\begin{lemm+}}{\end{lemm+}}

\usepackage{amsthm}
\theoremstyle{plain} \theoremstyle{remark}
\newtheorem{remark}{Remark}

\newcommand{\rank}{\mbox{rank}}
\def \r{\mbox{${\mathbb R}$}}
\def\E{/\kern-1.0em \equiv }

\title{Biharmonic maps in two dimensions }

\author{Ye-Lin Ou$^{*}$ and Sheng Lu }
\address{Department of
Mathematics,\newline\indent Texas A $\&$ M University-Commerce,
\newline\indent Commerce TX 75429,\newline\indent USA.\newline\indent
E-mail:yelin$\_$ou@tamu-commerce.edu \newline\indent \vskip0.1cm
School of Mathematics and Computer Science,\newline\indent Guangxi
University for Nationalities,
\newline\indent 188 East Daxue Road,\newline\indent Nanning, Guangxi 530006,\newline\indent P. R. China}
\thanks{ \indent* Supported by Texas A $\&$ M University-Commerce ``Faculty
Research Enhancement Project" (2010-11)}
\date{08/03/10}
\begin{document}
\title[Biharmonic maps in two dimensions]
{Biharmonic maps in two dimensions}


\maketitle
\section*{Abstract}
\begin{quote}
{\footnotesize Biharmonic maps between surfaces are studied in this
paper. We compute the bitension field of a map between surfaces with
conformal metrics in complex coordinates. As applications, we show
that a linear map from Euclidean plane into $(\mathbb{R}^2,
\sigma^2dwd\bar w)$ is always biharmonic if the conformal factor
$\sigma$ is bi-analytic; we construct a family of such $ \sigma$,
and we give a classification of linear biharmonic maps between
$2$-spheres. We also study biharmonic maps between surfaces with
warped product metrics. This includes a classification of linear
biharmonic maps between hyperbolic planes and some constructions of
many proper biharmonic maps into a circular cone or a helicoid.}
\end{quote}
\section{introduction and preliminaries}
All manifolds, maps, and tensor fields studied in this paper are assumed to be smooth.\\

A map $\phi : (M,g) \longrightarrow (N,h)$  between Riemannian
manifolds is {\em harmonic} if it is a critical point of the energy
functional
$$ E(\phi,\Omega ) = \frac{1}{2} \int_{M} |d\phi|^2 \, v_{g} $$
for every compact subset $\Omega$ of $M$. It is well known that
$\phi$ is harmonic if and only if its tension field $ \tau(\phi) =
{\rm Trace}_g\,{ \nabla \,d \phi }$ vanishes identically.\\

A map $\phi:(M, g)\longrightarrow (N, h)$ between Riemannian
manifolds is {\em biharmonic} if it is a critical point of the
bienergy functional
\begin{equation}\nonumber
E^{2}\left(\phi,\Omega \right)= \frac{1}{2} {\int}_{\Omega}
\left|\tau(\phi) \right|^{2}\, v_{g}
\end{equation}
for every compact subset $\Omega$ of $M$, where $\tau(\phi)={\rm
Trace}_{g}\nabla {\rm d} \phi$ is the tension field of $\phi$. The
Euler-Lagrange equation of this functional gives the biharmonic map
equation (\cite{Ji1})
\begin{equation}\label{BTF}
\tau^{2}(\phi):={\rm
Trace}_{g}(\nabla^{\phi}\nabla^{\phi}-\nabla^{\phi}_{\nabla^{M}})\tau(\phi)
- {\rm Trace}_{g} R^{N}({\rm d}\phi, \tau(\phi)){\rm d}\phi =0,
\end{equation}
where $R^{N}$ denotes the curvature operator of $(N, h)$ defined by
$$R^{N}(X,Y)Z=
[\nabla^{N}_{X},\nabla^{N}_{Y}]Z-\nabla^{N}_{[X,Y]}Z.$$

It follows from (\ref{BTF}) that a harmonic map is always biharmonic
so we use {\em proper biharmonic maps} to mean those biharmonic maps
which are not harmonic. For some recent study of constructions and
classifications of proper biharmonic maps see \cite{BK}, \cite{BFO},
 \cite{BMO1} ,  \cite{BMO2},  \cite{CMO1},  \cite{CMO2}, \cite{IIU},  \cite{LO},
\cite{MO},  \cite{On},  \cite{Ou1},  \cite{Ou2} , \cite{Ou3},
\cite{Ou4} ,  \cite{OT},  \cite{Oua},  \cite{WO} and the references
therein.\\

Harmonic maps between surfaces had been studied by many authors
during 1970' and 80' (see e.g., \cite{EL},  \cite{EW}, \cite{Jo},
\cite{Le} , \cite{SY}). It is well known that harmonic maps in two
dimensions have many special features and many of their important
properties are derived from the fact that the tension field of the
map
\begin{equation}\label{cxmap}
\phi: (M^2,g=\rho^2dzd{\bar z})\longrightarrow (N^2,
h=\sigma^2dwd{\bar w}),\;\;w=\phi(z)
\end{equation}

can be written in the concise formula
\begin{eqnarray}\label{tau}
\tau (\phi)=4\rho^{-2}[\phi_{z {\bar z}}+2(\ln \sigma)_{
w}\phi_z\phi_{{\bar z}}].
\end{eqnarray}

In this paper, we first prove that the bitension field of the map
defined in (\ref{cxmap}) can be written as
\begin{eqnarray}\notag\label{BTS}
\tau^2(\phi)&=&4\rho^{-2}\{\tau_{z\bar{z}}+2(\ln\sigma)_w[\tau_z\phi_{\bar{z}}
+\tau_{\bar{z}}\phi_z+\frac{1}{4}\rho^2(\tau)^2]+2\bar{\tau}(\ln\sigma)_{w\bar{w}}\phi_z\phi_{\bar{z}}\\
&&+2\tau\phi_z\phi_{\bar{z}}(\ln\sigma)_{ww}\},
\end{eqnarray}
where $\tau$ given by (\ref{tau}) is the tension field of $\phi$.\\

We then use this formula to show that a linear map from Euclidean
plane into $(\mathbb{R}^2, \sigma^2dwd\bar w)$ is always biharmonic
if the conformal factor $\sigma$ is bi-analytic (i.e.,
$\sigma_{ww}=0$); we construct a family of such $ \sigma$, and we
give a classification of linear biharmonic maps between $2$-spheres
using the model
$S^2\setminus\{N\}=(\mathbb{R}^2,\frac{4(dx^2+dy^2)}{(1+x^2+y^2)^2})$.
The second part of the paper is devoted to the study of  biharmonic
maps between surfaces with warped product metrics. Our results
include a classification of linear biharmonic maps between
hyperbolic planes and some constructions of many proper biharmonic
maps into a circular cone or a helicoid. The proper biharmonic maps
between hyperbolic planes found in our classification theorem  can
be viewed as a complement to a theorem proved in \cite{On} which
asserts that if $\phi : (M, g) \longrightarrow (N, h)$ is a map with
the property that $|\tau(\phi)| =$ constant, ${\rm Riem}^N < 0$, and
there is a point $p \in M$ such that $\rank_p\, \phi\ge 2$, then
$\phi$ is biharmonic if and only if it is harmonic.\\

As a preliminary step, we give a proof of the following lemma which
will be used frequently in the rest of the paper.
\begin{lemma}\label{OL}
Let $\phi :(M^{m}, g)\longrightarrow (N^{n}, h)$ be a map between
Riemannian manifolds with $\phi (x^{1},\ldots, x^{m})=(\phi^{1}(x),
\ldots, \phi^{n}(x))$ with respect to local coordinates $(x^{i})$ in
$M$ and $(y^{\alpha})$ in $N$. Then, $\phi$ is biharmonic if and
only if it is a solution of the following system of PDE's
\begin{eqnarray}\notag\label{BI3}
&&   \Delta\tau^{\sigma} +2g(\nabla\tau^{\alpha},
\nabla\phi^{\beta}) {\bar \Gamma_{\alpha\beta}^{\sigma}}
+\tau^{\alpha}\Delta \phi ^{\beta}{\bar
\Gamma_{\alpha\beta}^{\sigma}}\\ && +\tau^{\alpha}
g(\nabla\phi^{\beta}, \nabla\phi^{\rho})(\partial_{\rho}{\bar
\Gamma_{\alpha\beta}^{\sigma}}+{\bar
\Gamma_{\alpha\beta}^{\nu}}{\bar \Gamma_{\nu\rho}^{\sigma}})
-\tau^{\nu}g(\nabla\phi^{\alpha}, \nabla\phi^{\beta}){\bar
R}_{\beta\,\alpha \nu}^{\sigma}=0,\\\notag && \sigma=1,\, 2,\,
\ldots, n,
\end{eqnarray}
where $\tau^1, \ldots, \tau^n$ are components of the tension field
of the map $\phi$, $\nabla,\;\Delta$ denote the gradient and the
Laplace operators defined by the metric $g$, and
${\bar\Gamma_{\alpha\beta}^{\sigma}}$ and ${\bar R}_{\beta\,\alpha
\nu}^{\sigma}$ are the components of the connection and the
curvature of the target manifold.
\end{lemma}
\begin{proof}
Let $\{\frac{\partial}{\partial x^{i}}\}$ (resp.
$\{\frac{\partial}{\partial y^{\alpha}}\}$) be the natural frame
with respect to local coordinates $(x^{i})$ in $M$ (resp.
$(y^{\alpha})$ in $N$). Then, by (\ref{BTF}), the bitension field of
$\phi$ can be computed as

\begin{eqnarray}\label{TAU}
\tau^{2}(\phi) &=&g^{ij}\left(
\nabla^{\phi}_{\frac{\partial}{\partial
x^{i}}}\nabla^{\phi}_{\frac{\partial}{\partial x^{j}}} \left(
\tau(\phi)\right) -  \Gamma^{k}_{ij}
\nabla^{\phi}_{\frac{\partial}{\partial x^{k}}}
 \left( \tau(\phi)\right)\right) \\\notag&&- g^{ij}{\phi^{\alpha}}_{i} \phi^{\beta}_{j}
{\bar R}(\frac{\partial}{\partial y^{\alpha}},
\tau(\phi))\frac{\partial}{\partial y^{\beta}}.
\end{eqnarray}
Let $\tau(\phi)=\tau^{\alpha}\frac{\partial}{\partial y^{\alpha}}$
and use the notations $\tau^{\alpha}_i:= \frac{\partial
\tau^{\alpha}}{\partial x^{i}},\;\;\tau^{\alpha}_{ij}:=
\frac{\partial ^2\tau^{\alpha}}{\partial x^{i}\partial x^{j}}$.
Then, a direct computation gives
\begin{eqnarray}\label{T1}
\nabla^{\phi}_{\frac{\partial}{\partial x^{k}}}
 \left( \tau(\phi)\right)=\left(
\tau^{\sigma}_{k}+\tau^{\alpha}\phi^{\beta}_{k}{\bar
\Gamma_{\alpha\beta}^{\sigma}} \right)\frac{\partial}{\partial
y^{\sigma}},
\end{eqnarray}
\begin{eqnarray}\notag
\nabla^{\phi}_{\frac{\partial}{\partial
x^{i}}}\nabla^{\phi}_{\frac{\partial}{\partial x^{j}}} \left(
\tau(\phi)\right)&=&\nabla^{\phi}_{\frac{\partial}{\partial
x^{i}}}\left( \tau^{\alpha}_{j}+\tau^{\nu}\phi^{\beta}_{j}{\bar
\Gamma_{\beta\nu}^{\alpha}} \right)\frac{\partial}{\partial
y^{\alpha}}\\\label{T3} &=& \left( \tau^{\sigma}_{ij}
+\tau^{\alpha}_{j}\phi^{\beta}_{i}{\bar
\Gamma_{\alpha\beta}^{\sigma}} +\frac{\partial}{\partial
x^{i}}(\tau^{\alpha}\phi ^{\beta}_{j}{\bar
\Gamma_{\alpha\beta}^{\sigma}}) +\tau^{\alpha}\phi
^{\beta}_{j}\phi^{\rho}_{i}{\bar \Gamma_{\alpha\beta}^{\nu}}{\bar
\Gamma_{\nu\rho}^{\sigma}} \right)\frac{\partial}{\partial
y^{\sigma}},
\end{eqnarray}
and \begin{eqnarray}\label{T2} {\bar R}(\frac{\partial}{\partial
y^{\alpha}}, \tau(\phi))\frac{\partial}{\partial y^{\beta}}=
\tau^{\nu}{\bar R}_{\beta\,\alpha
\nu}^{\sigma}\frac{\partial}{\partial y^{\sigma}}.
\end{eqnarray}
 Substituting Equations (\ref{T1}), (\ref{T3}), and (\ref{T2}) into
(\ref{TAU}) we have

\begin{eqnarray}\notag
\tau^{2}(\phi) &=&   g^{ij}\{\tau^{\sigma}_{ij}
+\tau^{\alpha}_{j}\phi^{\beta}_{i}{\bar
\Gamma_{\alpha\beta}^{\sigma}} +\frac{\partial}{\partial
x^{i}}(\tau^{\alpha}\phi ^{\beta}_{j}{\bar
\Gamma_{\alpha\beta}^{\sigma}}) +\tau^{\alpha}\phi
^{\beta}_{j}\phi^{\rho}_{i}{\bar \Gamma_{\alpha\beta}^{\nu}}{\bar
\Gamma_{\nu\rho}^{\sigma}}\\
&&
-\Gamma^{k}_{ij}(\tau^{\sigma}_{k}+\tau^{\alpha}\phi^{\beta}_{k}{\bar
\Gamma_{\alpha\beta}^{\sigma}}) -\tau^{\nu}{\phi^{\alpha}}_{i}
\phi^{\beta}_{j}{\bar R}_{\beta\,\alpha
\nu}^{\sigma}\}\frac{\partial}{\partial y^{\sigma}},
\end{eqnarray}
from which the lemma follows.
\end{proof}

\section{Biharmonic maps between surfaces with conformal metrics}

Let $(M^2,g)$ and $(N^2, h)$ be two surfaces. As it is guaranteed by
a classical theorem from differential geometry we can always choose
local isothermal coordinates on a surface. So let $x,y$ on $M$
(respectively $u, v$ on $N$) be isothermal coordinates with respect
to which the metric takes the form $g=\rho^2\delta_{ij}dx^idx^j$
(respectively,
$h=\sigma^2\delta_{\alpha\beta}du^{\alpha}du^{\beta}$). Let $\phi:
(M^2,g)\longrightarrow (N^2, h)$,
$\phi(x,y)=(\phi^1(x,y),\phi^2(x,y))$ be a map between two surfaces.
It has been proved that in this case it is very useful to use
complex notations
\begin{eqnarray}\label{1}
&&z=x+iy, w=u+iv, \\\notag && dz=dx+idy,\;\; {\rm and}\;\;d{\bar
z}=dx-idy.\\\notag &&\frac{\partial}{\partial
z}=\frac{1}{2}(\frac{\partial}{\partial x}-i
\frac{\partial}{\partial y}),\;\; \frac{\partial}{\partial {\bar
z}}=\frac{1}{2}(\frac{\partial}{\partial x}+i
\frac{\partial}{\partial y}).
\end{eqnarray}
It is well known (see e.g., \cite{EL} and \cite{Jo}) that the
tension field of the map  $\phi: (M^2,g=\rho^2dzd{\bar
z})\longrightarrow (N^2, h=\sigma^2dwd{\bar w}),\;\; w=\phi(z)$ can
be written as

\begin{eqnarray}\label{CXtau}
\tau (\phi)=4\rho^{-2}[\phi_{z {\bar z}}+2(\ln \sigma)_{
w}\phi_z\phi_{{\bar z}}].
\end{eqnarray}

\begin{theorem}\label{MT}
Let $\phi: (M^2,g=\rho^2dzd{\bar z})\longrightarrow (N^2,
h=\sigma^2dwd{\bar w})$ with $w=\phi(z)$ be a map between two
surfaces with conformal metrics. Then, the bitension field of $\phi$
can be written as
\begin{eqnarray}\notag\label{BTS}
\tau^2(\phi)&=&4\rho^{-2}\{\tau_{z\bar{z}}+2(\ln\sigma)_w[\tau_z\phi_{\bar{z}}
+\tau_{\bar{z}}\phi_z+\frac{1}{4}\rho^2(\tau)^2]+2\bar{\tau}(\ln\sigma)_{w\bar{w}}\phi_z\phi_{\bar{z}}\\
&&+2\tau\phi_z\phi_{\bar{z}}(\ln\sigma)_{ww}\},
\end{eqnarray}
where $\tau$ given by (\ref{CXtau}) is the tension field of $\phi$.
\end{theorem}
\begin{proof}
A direct computation gives the connection coefficients of the target surface as\\
\begin{eqnarray}
&& \bar{\Gamma}^1_{11}=(\ln \sigma)_1,
\hskip0.3cm\bar{\Gamma}^1_{12}=(\ln \sigma)_2, \hskip0.3cm
\bar{\Gamma}^1_{22}=-(\ln \sigma)_1,\\\notag &&
\bar{\Gamma}^2_{11}=-(\ln \sigma)_2, \hskip0.3cm
\bar{\Gamma}^2_{12}=(\ln \sigma)_1, \hskip0.3cm
\bar{\Gamma}^2_{22}=(\ln \sigma)_2.\\\notag
\end{eqnarray}

Noting that the Laplace operator can be written as $\Delta
=4\rho^{-2}\frac{\partial^2}{\partial z\partial {\bar z}}$ we have
\begin{eqnarray}\label{G1}
&&\Delta\tau^1+i\Delta\tau^2=\Delta\tau=4\rho^{-2}\tau_{z\bar z}
\end{eqnarray}
Let $C_{\alpha\beta}=\begin{cases} 1,
\;\;\alpha=\beta=1\\i,\;\;\alpha\ne \beta\\-1, \;\;\alpha=\beta=2.
\end{cases}$
Then, it is easily checked that
\begin{equation}
\bar{\Gamma}^1_{\alpha\beta}
+i\bar{\Gamma}^2_{\alpha\beta}=2C_{\alpha\beta}(\ln \sigma)_w.
\end{equation}
We compute
\begin{eqnarray}\label{G2}
&&2g(\nabla\tau^{\alpha},\nabla\phi^{\beta})\bar{\Gamma}^1_{\alpha\beta}
+i2g(\nabla\tau^{\alpha},\nabla\phi^{\beta})\bar{\Gamma}^2_{\alpha\beta}\\\notag
=&&4(\ln
\sigma)_wg(\nabla\tau^{\alpha},\nabla\phi^{\beta})C_{\alpha\beta}=
8\rho^{-2}(\ln\sigma)_w(\tau_z\phi_{\bar{z}}+\tau_{\bar{z}}\phi_z),
\end{eqnarray}
\begin{eqnarray}\label{G3}
\tau^{\alpha}\Delta\phi^{\beta}\bar{\Gamma}^1_{\alpha\beta}
+i\tau^{\alpha}\Delta\phi^{\beta}\bar{\Gamma}^2_{\alpha\beta}
=2(\ln\sigma)_w\tau^{\alpha}\Delta\phi^{\beta}C_{\alpha\beta}
=8\rho^{-2}\tau(\ln\sigma)_w\phi_{z\bar z},
\end{eqnarray}
\begin{eqnarray}\label{G4}
&&\tau^{\alpha}g(\nabla\phi^{\beta},\nabla\phi^{\rho})\bar\Gamma^v_{\alpha\beta}\bar\Gamma^1_{v\rho}
+i\tau^{\alpha}g(\nabla\phi^{\beta},\nabla\phi^{\rho})\bar\Gamma^v_{\alpha\beta}\bar\Gamma^2_{v\rho}\\\notag
=&&4((\ln\sigma)_w)^2\tau\{g(\nabla\phi^{1},\nabla\phi^{1})+2ig(\nabla\phi^{1},\nabla\phi^{2})
-g(\nabla\phi^{2},\nabla\phi^{2})
\}\\\notag&&=16\tau\rho^{-2}((\ln\sigma)_w)^2\phi_z\phi_{\bar z},
\end{eqnarray}
and
\begin{eqnarray}
&&\tau^{\alpha}g(\nabla\phi^{\beta},\nabla\phi^{\rho})\partial_{\rho}\bar\Gamma^1_{\alpha\beta}
+i\tau^{\alpha}g(\nabla\phi^{\beta},\nabla\phi^{\rho})\partial_{\rho}\bar\Gamma^2_{\alpha\beta}
\\\notag
=&&2\tau\rho^2\{\partial_{1}((\ln\sigma)_w)[|\phi_{z}|^2+|\phi_{\bar
z}|^2)+2\phi_z\phi_{\bar
z}]\\\notag&&+\partial_{2}((\ln\sigma)_w)i[|\phi_{z}|^2+|\phi_{\bar
z}|^2)-2\phi_z\phi_{\bar z}]\}
\\\notag
=&&4\rho^{-2}\tau[(\ln\sigma)_{w\bar w}(|\phi_z|^2+|\phi_{\bar
z}|^2)+2(\ln\sigma)_{ww}\phi_z\phi_{\bar z}].
\end{eqnarray}

Using the fact that the Gauss curvature $K$ of $(N^2,
\sigma^2dwd\bar w)$ satisfies $K=-4\sigma^{-2}(\ln \sigma)_{w\bar
w}$ and hence $ \bar{R}^j_{i\;kl}=-4(\ln\sigma)_{w\bar
w}(\delta^j_k\delta_{il}-\delta^j_l\delta_{ik}) $ we have
\begin{eqnarray}\label{G6}
&&-\tau^vg(\nabla\phi^{\alpha},\nabla\phi^{\beta})\bar{R}^1_{\beta\alpha
v}-i\tau^vg(\nabla\phi^{\alpha},\nabla\phi^{\beta})\bar{R}^2_{\beta\alpha
v}\\\notag
=&&-\tau^1g(\nabla\phi^2,\nabla\phi^2)\bar{R}^1_{221}-\tau^2g(\nabla\phi^1,\nabla\phi^2)\bar{R}^1_{212}\\\notag
&&-i[\tau^2g(\nabla\phi^1,\nabla\phi^1)\bar{R}^2_{112}+\tau^1g(\nabla\phi^2,\nabla\phi^1)\bar{R}^2_{121}]\\\notag
=&&-4\rho^{-2}(\ln\sigma)_{w\bar w}[\tau(|\phi_z|^2+|\phi_{\bar
z}|^2)-2\bar{\tau}\phi_z\phi_{\bar z}].
\end{eqnarray}
Substituting Equations (\ref{G1}) and (\ref{G2})--(\ref{G6}) into
the bitension field equation  we obtain
\begin{eqnarray}\notag
&&\tau^2(\phi)=(\tau^2(\phi))^1+i(\tau^2(\phi))^2\\\notag
=&&\Delta\tau^1+2g(\nabla\tau^{\alpha},\nabla\phi^{\beta})\bar\Gamma^1_{\alpha\beta}+\tau^{\alpha}\Delta\phi^{\beta}\bar\Gamma^1_{\alpha\beta}\\\notag
&&+\tau^{\alpha}g(\nabla\phi^{\beta},\nabla\phi^{\rho})(\partial_{\rho}\bar\Gamma^1_{\alpha\beta}+\bar\Gamma^v_{\alpha\beta}\bar\Gamma^1_{v\rho})
-\tau^vg(\nabla\phi^{\alpha},\nabla\phi^{\beta})\bar{R}^1_{\beta\alpha
v}\\\notag
&&+i\,[\Delta\tau^2+2g(\nabla\tau^{\alpha},\nabla\phi^{\beta})\bar\Gamma^2_{\alpha\beta}
+\tau^{\alpha}\Delta\phi^{\beta}\bar\Gamma^2_{\alpha\beta}\\\notag
&&+\tau^{\alpha}g(\nabla\phi^{\beta},\nabla\phi^{\rho})(\partial_{\rho}\bar\Gamma^2_{\alpha\beta}
+\bar\Gamma^v_{\alpha\beta}\bar\Gamma^2_{v\rho})-\tau^vg(\nabla\phi^{\alpha},\nabla\phi^{\beta})\bar{R}^2_{\beta\alpha
v})]\\\notag
=&&4\rho^{-2}\{\tau_{z\bar{z}}+2(\ln\sigma)_w[\tau_z\phi_{\bar{z}}
+\tau_{\bar{z}}\phi_z+\frac{1}{4}\rho^2(\tau)^2]+2\bar{\tau}(\ln\sigma)_{w\bar{w}}\phi_z\phi_{\bar{z}}\\\notag
&&+2\tau(\ln\sigma)_{ww}\phi_z\phi_{\bar{z}}\}.
\end{eqnarray}
From this the theorem follows.
\end{proof}
\begin{corollary}\label{BLE}
A map $\phi: (M^2,g=\rho^2dzd{\bar z})\longrightarrow (N^2,
h=\sigma^2dwd{\bar w})$ with $w=\phi(z)$ between surfaces is
biharmonic if and only if
\begin{eqnarray}\label{BSE}
&&\tau_{z\bar{z}}+2(\ln\sigma)_w[\tau_z\phi_{\bar{z}}
+\tau_{\bar{z}}\phi_z+\frac{1}{4}\rho^2(\tau)^2]\\\notag
&&+2\phi_z\phi_{\bar{z}}[\bar{\tau}(\ln\sigma)_{w\bar{w}}
+\tau(\ln\sigma)_{ww}]=0,
\end{eqnarray}
where $\tau (\phi)=4\rho^{-2}[\phi_{z {\bar z}}+2(\ln \sigma)_{
w}\phi_z\phi_{{\bar z}}]$ is the tension field of the map $\phi$.
\end{corollary}
As an application of Theorem \ref{MT}, we prove the following
theorem that gives a lot of examples of linear maps which are proper
biharmonic maps.
\begin{theorem}\label{P0}
Let $\phi:(\mathbb{R}^2,dx^2+dy^2)\longrightarrow
(\mathbb{R}^2,\sigma^2(du^2+dv^2))$ with
\begin{equation}\label{MAP}
\phi(x,y)=\left(\begin{array}{cc} a & b\\c& d
\end{array}\right)\left(\begin{array}{c} x\\y
\end{array}\right)
\end{equation}
 be a linear map. Then \\ (1)
$\phi$ is harmonic if and only if $\sigma$ is anti-analytic
($\sigma_w=0$), or $\phi$ is
analytic ($\phi_{\bar z}=0)$  or anti-analytic ($\phi_z=0$);\\
(2) If $\sigma$ is anti-bianalytic (i.e., $\sigma_{ww}=0$), then the
linear map $\phi$ is always a biharmonic map. In particular, for any
$p,q \in \r, q\ne 0$, the linear map
$\phi:(\mathbb{R}^2,dx^2+dy^2)\longrightarrow
(\mathbb{R}^2,(p+q(u^2+v^2))^2(du^2+dv^2))$ defined by (\ref{MAP})
with $|A^1|^2\neq |A^2|^2$ or $A^1\cdot A^2\neq0$, is always a
proper biharmonic map, where $A^1=(a,b)$ and $A^2=(c,d)$.
\end{theorem}
\begin{proof}
The map can be written as\\
\begin{equation}\label{Lmap}
\phi(z)=\frac{1}{2}(a+ic)(z+{\bar z})+\frac{1}{2}(d-ib)(z-{\bar z}),
\end{equation}
which is  linear in $z$ and $\bar z$. With this, together with the
fact that $\rho=1$, we can compute the tension field of the map to
get
\begin{eqnarray}\label{abc}
&&\tau=8(\ln\sigma)_w\phi_z\phi_{\bar z},
\end{eqnarray}
from which we obtain the first statement of the theorem.\\

To prove the second statement of the theorem we compute
\begin{eqnarray}\label{30}
\tau_{z\bar z}=&&8\phi_z\phi_{\bar
z}((\ln\sigma)_{www}\phi_z\phi_{\bar z}+(\ln\sigma)_{ww\bar
w}(\phi_z\bar{\phi}_{\bar z}+\phi_{\bar z}\bar{\phi}_z)\\\notag
&&+(\ln\sigma)_{w\bar w\bar w}\bar{\phi}_z\bar{\phi}_{\bar z}),
\end{eqnarray}
\begin{eqnarray}
&&2(\ln\sigma)_w(\tau_z\phi_{\bar z}+\tau_{\bar
z}\phi_z+\frac{1}{4}\rho^2(\tau)^2)\\\notag
&&=16(\ln\sigma)_w\phi_z\phi_{\bar
z}\big[2(\ln\sigma)_{ww}\phi_z\phi_{\bar z}+(\ln\sigma)_{w\bar
w}(\phi_{\bar z}\bar{\phi}_z+\phi_z\bar{\phi}_{\bar
z})+2(\ln\sigma)^2_w\phi_z\phi_{\bar z}\big],
\end{eqnarray}
\begin{eqnarray}
2(\ln\sigma)_{w\bar w}\bar{\tau}\phi_z\phi_{\bar z}
&=&16(\ln\sigma)_{\bar w}(\ln\sigma)_{w\bar
w}\bar{\phi}_z\bar{\phi}_{\bar z}\phi_z\phi_{\bar z},
\end{eqnarray}
and
\begin{eqnarray}\label{31}
2\tau\phi_z\phi_{\bar
z}(\ln\sigma)_{ww}=16(\ln\sigma)_w(\ln\sigma)_{ww}\phi_z^2\phi_{\bar
z}^2.
\end{eqnarray}
Substituting Equations (\ref{30})$\sim$(\ref{31}) into biharmonic
map equation (\ref{BSE}) we see that $\phi$ is biharmonic if and
only if
\begin{eqnarray}\notag
&&\phi_z\phi_{\bar z}[8(\ln\sigma)_{www}\phi_z\phi_{\bar
z}+8(\ln\sigma)_{ww\bar w}(\phi_z\bar{\phi}_{\bar z}+\phi_{\bar
z}\bar{\phi}_z)+8(\ln\sigma)_{w\bar w\bar
w}\bar{\phi}_z\bar{\phi}_{\bar z}\\\label{32}
&&+48(\ln\sigma)_w(\ln\sigma)_{ww}\phi_z\phi_{\bar
z}+16(\ln\sigma)_w(\ln\sigma)_{w\bar w}(\phi_{\bar
z}\bar{\phi}_z+\phi_z\bar{\phi}_{\bar z})\\\notag
&&+16(\ln\sigma)_{\bar w}(\ln\sigma)_{w\bar
w}\bar{\phi}_z\bar{\phi}_{\bar
z}+32((\ln\sigma)_w)^3\phi_z\phi_{\bar z}]=0.
\end{eqnarray}
Noting that $\phi_z\phi_{\bar z}=0$ implies $\phi$ is harmonic we
conclude that if $\phi$ is proper biharmonic, then, by (\ref{32}),
it solves the equation
\begin{eqnarray}\notag
&&8(\ln\sigma)_{www}\phi_z\phi_{\bar z}+8(\ln\sigma)_{ww\bar
w}(\phi_z\bar{\phi}_{\bar z}+\phi_{\bar
z}\bar{\phi}_z)+8(\ln\sigma)_{w\bar w\bar
w}\bar{\phi}_z\bar{\phi}_{\bar z}\\\label{33}
&&+48(\ln\sigma)_w(\ln\sigma)_{ww}\phi_z\phi_{\bar
z}+16(\ln\sigma)_w(\ln\sigma)_{w\bar w}(\phi_{\bar
z}\bar{\phi}_z+\phi_z\bar{\phi}_{\bar z})\\\notag
&&+16(\ln\sigma)_{\bar w}(\ln\sigma)_{w\bar
w}\bar{\phi}_z\bar{\phi}_{\bar
z}+32((\ln\sigma)_w)^3\phi_z\phi_{\bar z}=0.
\end{eqnarray}
It is not difficult to check that Equation (\ref{33}) is equivalent
to
\begin{eqnarray}\notag
&&8\phi_z\phi_{\bar
z}[((\ln\sigma)_{ww}+((\ln\sigma)_w)^2)_w+4(\ln\sigma)_w((\ln\sigma)_{ww}+((\ln\sigma)_w)^2)]\\\label{BFE}
&&+8\bar{\phi}_z\bar{\phi}_{\bar
z}(\overline{(\ln\sigma)_{ww}+((\ln\sigma)_w)^2})_w\\\notag
&&+8(\phi_z\bar{\phi}_{\bar z}+\phi_{\bar
z}\bar{\phi}_z)((\ln\sigma)_{ww}+((\ln\sigma)_w)^2)_{\bar w}=0.
\end{eqnarray}
It follows that if
\begin{eqnarray}\label{34}
(\ln\sigma)_{ww}+((\ln\sigma)_w)^2=0
\end{eqnarray}
then  (\ref{BFE}), and hence (\ref{33}), is
automatically solved.\\

Since Equation (\ref{34}) is equivalent to
\begin{eqnarray}\label{35}
&&\sigma_{ww}=0,
\end{eqnarray}
from which we obtain the first part of Statement (2). It is easy to
check that $\sigma=p+qw\bar w=p+q(u^2+v^2)$ is a special solution of
(\ref{35}). It follows that, any linear map
$\phi:(\mathbb{R}^2,dx^2+dy^2)\longrightarrow
(\mathbb{R}^2,(p+q(u^2+v^2))^2(du^2+dv^2))$ with
$\phi(x,y)=(ax+by,cx+dy)$ is biharmonic map; Furthermore, if
$|A^1|^2\neq |A^2|^2$ or $A^1\cdot A^2\neq0$, and hence
$\phi_z\phi_{\bar z}=\frac{1}{4}(|A^1|^2-|A^2|^2+i2A^1\cdot A^2)\ne
0$ which means $\phi$ is not harmonic, then $\phi$ is a proper
biharmonic map. This completes the proof of the theorem.
\end{proof}
As another application of Theorem \ref{MT}, we give a classification
of linear biharmonic maps between $2$-spheres with the model
$S^2\setminus\{N\}=(\mathbb{R}^2,\frac{4(dx^2+dy^2)}{(1+x^2+y^2)^2})$.
\begin{proposition}\label{P0}
A linear map between $2$-spheres
$\phi:S^2\setminus\{N\}=(\mathbb{R}^2,\frac{4(dx^2+dy^2)}{(1+x^2+y^2)^2})\longrightarrow
S^2\setminus\{N\}=(\mathbb{R}^2,\frac{4(du^2+dv^2)}{(1+u^2+v^2)^2})$
with $\phi(x,y)=(ax+by,cx+dy)$ is biharmonic if and only if it is
harmonic, i.e., a constant map, or $|A^1|^2=|A^2|^2$ and $A^1\cdot
A^2=0$, where $A^1=(a,b), \;\; A^2=(c,d)$ .
\end{proposition}
\begin{proof}
The linear map can be put in the form of (\ref{Lmap}). Here, we have
\begin{eqnarray*}
&&\rho=\frac{2}{1+z\bar z},\;\; \sigma=\frac{2}{1+w\bar w},
\end{eqnarray*}
and the tensions field of $\phi$ given by
\begin{eqnarray}\notag
\tau (\phi)=4\rho^{-2}[\phi_{z {\bar
z}}+2(\ln\sigma)_w\phi_z\phi_{{\bar z}}]=8\rho^{-2}(\ln \sigma)_{
w}\phi_z\phi_{{\bar z}}.
\end{eqnarray}

A straightforward computation gives
\begin{eqnarray}\label{22}
&&\tau_{z\bar z}=4[\rho^{-2}(\phi_{z\bar
z}+2(\ln\sigma)_w\phi_z\phi_{\bar z})]_{z\bar z}\\\notag
=&&\frac{8}{(1+w\bar w)^3}[-\frac{3}{2}z\bar z\bar w(1+w\bar
w)^2\phi_z\phi_{\bar z}-\frac{1-z\bar z}{2}\bar w(1+w\bar
w)^2\phi_z\phi_{\bar z}\\\notag &&+\frac{(1+z\bar z)\bar z}{2}\bar
w^2(1+w\bar w)\phi_z\phi_{\bar z}^2-\frac{(1+z\bar z)\bar
z}{2}(1+w\bar w)\phi_z\phi_{\bar z}\bar{\phi}_{\bar z}\\\notag
&&+\frac{(1+z\bar z)z}{2}\bar w^2(1+w\bar w)\phi_z^2\phi_{\bar
z}-\frac{(1+z\bar z)^2}{2}\bar w^3\phi_z^2\phi_{\bar
z}^2+\frac{(1+z\bar z)^2}{2}\bar w\phi_z^2\phi_{\bar
z}\bar{\phi}_{\bar z}\\\notag &&-\frac{(1+z\bar z)z}{2}(1+w\bar
w)\phi_z\phi_{\bar z}\bar{\phi}_z+\frac{(1+z\bar z)^2}{2}\bar
w\phi_z\phi_{\bar z}^2\bar{\phi}_z+\frac{(1+z\bar
z)^2}{2}w\phi_z\phi_{\bar z}\bar{\phi}_z\bar{\phi}_{\bar z}],
\end{eqnarray}

\begin{eqnarray}
&&2(\ln\sigma)_w(\tau_z\phi_{\bar z}+\tau_{\bar
z}\phi_z+\frac{1}{4}\rho^2(\tau)^2)\\\notag =&&2\frac{-\bar
w}{1+w\bar w}[(8\rho^{-2}(\ln\sigma)_w\phi_z\phi_{\bar
z})_z\phi_{\bar z}+(8\rho^{-2}(\ln\sigma)_w\phi_z\phi_{\bar
z})_{\bar z}\phi_z\\\notag
&&+16\rho^{-2}((\ln\sigma)_w)^2\phi_z^2\phi_{\bar z}^2]\\\notag =
&&-16\frac{1}{(1+w\bar w)^3}[-\frac{\bar z(1+z\bar z)}{2}\bar
w^2(1+w\bar w)\phi_z\phi_{\bar z}^2+\frac{(1+z\bar z)^2}{2}\bar
w^3\phi_z^2\phi_{\bar z}^2\\\notag &&-\frac{(1+z\bar z)^2}{4}\bar
w\phi_z\phi_{\bar z}^2\bar{\phi}_z-\frac{z(1+z\bar z)}{2}\bar
w^2(1+w\bar w)\phi_z^2\phi_{\bar z}-\frac{(1+z\bar z)^2}{4}\bar
w\phi_z^2\phi_{\bar z}\bar{\phi}_{\bar z}\\\notag &&+\frac{(1+z\bar
z)^2}{2}\bar w^3\phi_z^2\phi_{\bar z}^2],
\end{eqnarray}
\begin{eqnarray}
2(\ln\sigma)_{w\bar w}\bar{\tau}\phi_z\phi_{\bar
z}=\frac{4}{(1+w\bar w)^3}(1+z\bar z)^2w\phi_z\phi_{\bar
z}\bar{\phi}_z\bar{\phi}_{\bar z},\\\notag
\end{eqnarray}
and
\begin{eqnarray}\label{23}
2\tau\phi_z\phi_{\bar z}(\ln\sigma)_{ww}=4(1+z\bar z)^2\frac{-\bar
w^3}{(1+w\bar w)^3}\phi_z^2\phi_{\bar z}^2.\\\notag
\end{eqnarray}
Summing up (\ref{22}) $\sim$ (\ref{23}) and using Corollary
\ref{BLE} we conclude that $\phi$ is biharmonic if and only if
\begin{eqnarray}\label{24}
&&\frac{1}{(1+w\bar w)^3}[-12z\bar z\bar w(1+w\bar
w)^2\phi_z\phi_{\bar z}-4(1-z\bar z)\bar w(1+w\bar
w)^2\phi_z\phi_{\bar z}\\\notag &&+4(1+z\bar z)\bar{z}\bar
w^2(1+w\bar w)\phi_z\phi_{\bar z}^2-4(1+z\bar z)\bar{z}(1+w\bar
w)\phi_z\phi_{\bar z}\bar{\phi}_{\bar z}\\\notag &&+4(1+z\bar
z)z\bar w^2(1+w\bar w)\phi_z^2\phi_{\bar z}\\\notag &&-4(1+z\bar
z)^2\bar w^3\phi_z^2\phi_{\bar z}^2+4(1+z\bar z)^2\bar
w\phi_z^2\phi_{\bar z}\bar{\phi}_{\bar z}-4(1+z\bar z)z(1+w\bar
w)\phi_z\phi_{\bar z}\bar{\phi}_z\\\notag &&+4(1+z\bar z)^2\bar
w\phi_z\phi_{\bar z}^2\bar{\phi}_z+4(1+z\bar z)^2w\phi_z\phi_{\bar
z}\bar{\phi}_z\bar{\phi}_{\bar z}\\\notag &&+8\bar z(1+z\bar z)\bar
w^2(1+w\bar w)\phi_z\phi_{\bar z}^2-8(1+z\bar z)^2\bar
w^3\phi_z^2\phi_{\bar z}^2\\\notag &&+4(1+z\bar z)^2\bar
w\phi_z\phi_{\bar z}^2\bar{\phi}_z+8z(1+z\bar z)\bar w^2(1+w\bar
w)\phi_z^2\phi_{\bar z}+4(1+z\bar z)^2\bar w\phi_z^2\phi_{\bar
z}\bar{\phi}_{\bar z}\\\notag &&+4(1+z\bar z)^2w\phi_z\phi_{\bar
z}\bar{\phi}_z\bar{\phi}_{\bar z}-12(1+z\bar z)^2\bar
w^3\phi_z^2\phi_{\bar z}^2]=0.
\end{eqnarray}
Noting that $w=\phi(z)$ is linear in $z$ and $\bar z$ (and hence
$\phi_z$ and $\phi_{\bar z}$ are constants) we multiply $(1+w\bar
w)^3$ to both sides of (\ref{24}) to obtain a polynomial equation in
$ x, y\;\; (z=x+iy)$ whose 7th degree terms gives
\begin{eqnarray}
&&-12z\bar zw^2\bar w^3\phi_z\phi_{\bar z}+4z\bar zw^2\bar
w^3\phi_z\phi_{\bar z}+4z\bar z^2w\bar w^3\phi_z\phi_{\bar
z}^2\\\notag &&+4z^2\bar zw\bar w^3\phi_z^2\phi_{\bar z}-4z^2\bar
z^2\bar w^3\phi_z^2\phi_{\bar z}^2\\\notag &&+8z\bar z^2w\bar
w^3\phi_z\phi_{\bar z}^2-8z^2\bar z^2\bar w^3\phi_z^2\phi_{\bar
z}^2+8z^2\bar zw\bar w^3\phi_z^2\phi_{\bar z}\\\notag &&-12z^2\bar
z^2\bar w^3\phi_z^2\phi_{\bar z}^2=0.
\end{eqnarray}
This, by a straightforward computation, is equivalent to

\begin{eqnarray}
4\phi_z\phi_{\bar z}z\bar z\bar w^3(-2w^2+3\bar zw\phi_{\bar
z}+3zw\phi_z-6z\bar z\phi_z\phi_{\bar z})=0
\end{eqnarray}
for all $z\in \mathbb{C}$. It follows that either (i)
$w=\phi(z)\equiv 0$ which means $\phi$ is a constant map, or (ii)
$\phi_z\phi_{\bar z}=0$, which means the map $\phi$ is harmonic, or
(iii)
\begin{equation}\label{gu8}
-2w^2+3\bar zw\phi_{\bar z}+3zw\phi_z-6z\bar z\phi_z\phi_{\bar z}=0
\end{equation}
for all $z\in \mathbb{C}$.\\

Using $z=x+iy$, $w=\phi(z)=ax+by+i(cx+dy)$ we can easily check that
Equation (\ref{gu8}) is equivalent to
\begin{eqnarray*}
&&-\frac{a^2}{2}x^2-\frac{3b^2}{2}x^2+\frac{c^2}{2}x^2+\frac{3d^2}{2}x^2+2abxy-2cdxy\\\notag
&&-\frac{3a^2}{2}y^2-\frac{b^2}{2}y^2+\frac{3c^2}{2}y^2+\frac{d^2}{2}y^2\\\notag
&&+i[-acx^2-3bdx^2+2bcxy+2adxy-3acy^2-bdy^2]=0,
\end{eqnarray*}
or
\begin{equation}\label{25}
\begin{cases}
-\frac{1}{2}a^2-\frac{3}{2}b^2+\frac{1}{2}c^2+\frac{3}{2}d^2=0,\\
\\
2ab-2cd=0,\\
\\
-\frac{3}{2}a^2-\frac{1}{2}b^2+\frac{3}{2}c^2+\frac{1}{2}d^2=0,\\
\\
-ac-3bd=0,\\
\\
2bc+2ad=0,\\
\\
-3ac-bd=0.
\end{cases}
\end{equation}
Solving Equation (\ref{25}) we see that $a=b=c=d=0$ is the only
solution. Summarizing the above results we obtain the proposition.
\end{proof}

\section{ Biharmonic maps between surfaces with warped product metrics}
In this section, we study biharmonic maps from a Euclidean plane
into a surface with a warped product metric following the idea of
Lemaire \cite{Le} in his search for harmonic map $\phi:
T^2\longrightarrow (\mathbb{R}^2, du^2+(a^2-u^2)dv^2)$ of the form
$\phi(x,y)=(f(x),y)$. It turns out that  the partial differential
equations of a biharmonic map in this case reduce to an ordinary
differential equation. The solutions of the resulting ordinary
differential equation allow us to construct many proper biharmonic
maps from a Euclidean plane into a circular cone and a helicoid. We
also give a complete classification of linear biharmonic maps
between hyperbolic planes with the model
$H^2=(\mathbb{R}^2,e^{-2y}dx^2+dy^2)$.

\begin{proposition}\label{P0}
Let $\mathbb{R}^2$ be the Euclidean plane with coordinates $(x,y)$
and metric $g=dx^2+dy^2$. Let $N =\mathbb{R}^2$ with coordinates
$(u,v)$ be provided with the warped product metric
$d\rho^2=du^2+\sigma^2(u)dv^2$. Then, the map $\phi:
\mathbb{R}^2\longrightarrow N$ defined by $\phi(x,y)=(f(x),y)$
 is biharmonic if and only if
\begin{equation}\label{SE1}
(f''-\sigma'(f)\sigma(f))''-(f''-\sigma'(f)\sigma(f))(\sigma'\sigma)'(f)=0.
\end{equation}
\end{proposition}
\begin{proof}

For the target surface, we have the coefficients of the first
fundamental form $E=1, F=0, G=\sigma^2(u)$ and hence the Christoffel
symbols given by
\begin{eqnarray}
&&\bar{\Gamma}^1_{11}=E_u/2E=0,
\hskip0.3cm\bar{\Gamma}^1_{12}=E_v/2E=0, \hskip0.3cm
\bar{\Gamma}^1_{22}=-G_u/2E=-\sigma'\sigma,\\\notag &&
\bar{\Gamma}^2_{11}=-E_v/2G=0, \hskip0.3cm
\bar{\Gamma}^2_{12}=G_u/2G=\sigma'/\sigma, \hskip0.3cm
\bar{\Gamma}^2_{22}=G_v/2G=0.\\\notag
\end{eqnarray}
We can further compute the components of the tension field of the
map $\phi$ to get
\begin{eqnarray}
&&\tau^1=f''-\sigma'(f)\sigma(f),\\
&&\tau^2=(f')^2\bar{\Gamma}^2_{11}+\bar{\Gamma}^2_{22}=0.
\end{eqnarray}
A straightforward computation using Lemma \ref{OL} gives the
components of the bitension field of $\phi$ as:
\begin{eqnarray}\label{CPL1}
&&(\tau^2(\phi))^1= \Delta\tau^{1} +2g(\nabla\tau^{\alpha},
\nabla\phi^{\beta}) {\bar \Gamma_{\alpha\beta}^{1}}
+\tau^{\alpha}\Delta \phi ^{\beta}{\bar
\Gamma_{\alpha\beta}^{1}}\\\notag && +\tau^{\alpha}
g(\nabla\phi^{\beta}, \nabla\phi^{\rho})(\partial_{\rho}{\bar
\Gamma_{\alpha\beta}^{1}}+{\bar \Gamma_{\alpha\beta}^{\nu}}{\bar
\Gamma_{\nu\rho}^{1}}) -\tau^{\nu}g(\nabla\phi^{\alpha},
\nabla\phi^{\beta}){\bar R}_{\beta\,\alpha \nu}^{1}\\\notag =&&
\tau^{1}_{11} +\tau^{1}{\bar \Gamma_{12}^{2}}{\bar \Gamma_{22}^{1}}
-\tau^{1}(-\partial_1{\bar
\Gamma_{22}^{1}}+{\bar\Gamma_{12}^{2}}{\bar
\Gamma_{22}^{1}})\\\notag =&&
(f''-\sigma'(f)\sigma(f))''-(f''-\sigma'(f)\sigma(f))(\sigma'\sigma)'(f),
\end{eqnarray}
\begin{eqnarray}\label{CPL2}
&&(\tau^2(\phi))^2= \Delta\tau^{2} +2g(\nabla\tau^{\alpha},
\nabla\phi^{\beta}) {\bar \Gamma_{\alpha\beta}^{2}}
+\tau^{\alpha}\Delta \phi ^{\beta}{\bar
\Gamma_{\alpha\beta}^{2}}\\\notag && +\tau^{\alpha}
g(\nabla\phi^{\beta}, \nabla\phi^{\rho})(\partial_{\rho}{\bar
\Gamma_{\alpha\beta}^{2}}+{\bar \Gamma_{\alpha\beta}^{\nu}}{\bar
\Gamma_{\nu\rho}^{2}}) -\tau^{\nu}g(\nabla\phi^{\alpha},
\nabla\phi^{\beta}){\bar R}_{\beta\,\alpha \nu}^{2}\\\notag =&&
2g(\nabla\tau^{1}, \nabla\phi^{\beta}) {\bar \Gamma_{1\beta}^{2}}
+\tau^{1}\Delta \phi ^{\beta}{\bar \Gamma_{1\beta}^{2}} +\tau^{1}
g(\nabla\phi^{\beta}, \nabla\phi^{\rho})(\partial_{\rho}{\bar
\Gamma_{1\beta}^{2}}\\\notag &&+{\bar \Gamma_{1\beta}^{\nu}}{\bar
\Gamma_{\nu\rho}^{2}}) -\tau^{1}g(\nabla\phi^{\alpha},
\nabla\phi^{\beta}){\bar R}_{\beta\,\alpha 1}^{2}=0.
\end{eqnarray}
It follows that $\phi$ is biharmonic if and only if $f$ solves the
ordinary differential equation (\ref{SE1}), which completes the
proof of the proposition.
\end{proof}

\begin{corollary}\label{P1}
Let $a>0$ and $A, C$ be constants satisfying $|A|<a, C\ne 0$ and $B,
D$ be arbitrary constants. Let $M=\{(x,y)\in \mathbb{R}^2:
|x|<\frac{a-|A|}{|C|}\}$ be provided with the standard Euclidean
metric $ds^2=dx^2+dy^2$ and $N =\{(u,v)\in\mathbb{R}^2: |u|<a \}$ be
a surface with the  metric $d\rho^2=du^2+(a^2-u^2)dv^2$. Then, the
map\\ $\phi: (M,\;\; ds^2=dx^2+dy^2)\longrightarrow
(N,\;\;d\rho^2=du^2+(a^2-u^2)dv^2)$ defined by $\phi(x,y)=(A\cos
(x+B)+Cx\cos (x+D),\;\;y)$ is a proper biharmonic map.
\end{corollary}
\begin{proof}

To the effect, the map $\phi: M\longrightarrow N$ defined by
$\phi(x,y)=(f(x),y)$ can be viewed as a map obtained from the map
defined in Proposition \ref{P0} with $\sigma^2 (u)=a^2-u^2$ by
restricting its domain. It follows that $\sigma(u)\sigma'(u)=-u$,
and $(\sigma(u)\sigma'(u))'=-1$. Substituting these into  Equation
(\ref{SE1}) we see that the map $\phi$ is biharmonic if and only if
$f^{(iv)}+2f''+f=0$. Solving this 4th order linear equation with
constant coefficients we get $f(x)=A\cos (x+B)+Cx\cos (x+D)$. From
this we obtain the corollary.
\end{proof}
\begin{remark}
The idea in Corollary \ref{P1} was motivated by \cite{Le} in which
Lemaire used the equation $\tau
(\phi)=[f''-\sigma'(f)\sigma(f)]\frac{\partial}{\partial x}=0$ to
determine $f(x)=A\cos (x+B)$ which gives a continuous family of
harmonic maps which are not isometrically or conformally equivalent.
\end{remark}
\begin{corollary}\label{Helicoid}{\rm (Biharmonic maps into a
helicoid)} Let $\mathbb{R}^2$ be the Euclidean space with the
standard coordinates $(x,y)$ and the metric $ds^2=dx^2+dy^2$. Let
$N$ be the helicoid with the parametrization $\vec{r}(u,v)=(u\cos v,
u \sin v, av)$ and the induced metric from $\mathbb{R}^3$. Then, the
map $\phi:\mathbb{R}^2\longrightarrow N$ defined by
$\phi(x,y)=(f(x),y)$
 is biharmonic if and only if $f(x)=(A+Bx)e^x+(C+Dx)e^{-x}$. It is proper biharmonic when
$B^2+D^2\ne0$.
\end{corollary}
\begin{proof}
For helicoid $\vec{r}(u,v)=(u\cos v, u \sin v, av)$ we have $E=1,
F=0, G=u^2+a^2$, so the induced metric on the helicoid takes the
form $d\rho^2=du^2+(a^2+u^2)dv^2$. Therefore, the map $\phi:
\mathbb{R}^2\longrightarrow N$ defined by $\phi(x,y)=(f(x),y)$, or
$\begin{cases} u(x,y)=f(x) \\v(x,y)=y\end{cases}$, can be viewed as
a map obtained from the map defined in Proposition \ref{P0} with
$\sigma^2 (u)=a^2+u^2$. In this case, $\sigma(u)\sigma'(u)=u$, and
$(\sigma(u)\sigma'(u))'=1$. It follows from Equation (\ref{SE1})
that the map $\phi$ is biharmonic if and only if
$f^{(iv)}-2f''+f=0$. Solving this 4th order linear equation with
constant coefficients we get $f(x)=(A+Bx)e^x+(C+Dx)e^{-x}$.
Therefore, for $B^2+D^2\ne 0$, we obtain a family of proper
biharmonic maps from Euclidean plane into a helicoid.
\end{proof}

\begin{corollary}\label{Cone}{\rm (Biharmonic maps into a
circular cone)} Let $\mathbb{R}^2$ be the Euclidean space with the
standard coordinates $(x,y)$ and the metric $ds^2=dx^2+dy^2$. Let
$N$ be the circular cone in $\mathbb{R}^3$ with parametrization
$\vec{r}(u,v)=(\frac{u}{\sqrt{2}}\cos v, \frac{u}{\sqrt{2}} \sin v,
\frac{u}{\sqrt{2}})$ with the induced metric from $\mathbb{R}^3$.
Then, the map $\phi: T^2\longrightarrow N$ defined by
$\phi(x,y)=(f(x),y)$ is biharmonic if and only if
$f(x)=(A+Bx)e^{x/\sqrt{2}}+(C+Dx)e^{-x/\sqrt{2}}$. It is proper
biharmonic when $B^2+D^2\ne 0$.
\end{corollary}
\begin{proof}
For the circular cone in $\mathbb{R}^3$ with parametrization
$\vec{r}(u,v)=(\frac{u}{\sqrt{2}}\cos v, \frac{u}{\sqrt{2}} \sin v,
\frac{u}{\sqrt{2}})$, we have $E=1, F=0, G=\frac{1}{2}u^2$, so the
induced metric on the circular cone takes the form
$d\rho^2=du^2+\frac{1}{2}u^2dv^2$. Therefore, the map $\phi:
\mathbb{R}^2\longrightarrow N$ defined by $\phi(x,y)=(f(x),y)$ can
be viewed as a map obtained from the map defined in Proposition
\ref{P0} with $\sigma^2 (u)=\frac{1}{2}u^2$ and a restriction of the
co-domain. In this case, $\sigma(u)\sigma'(u)=\frac{u}{2}$, and
$(\sigma(u)\sigma'(u))'=1/2$, and it follows from Equation
(\ref{SE1}) that the map $\phi$ is biharmonic if and only if
$f^{(iv)}-f''+\frac{1}{4}f=0$. Solving this 4th order linear
equation with constant coefficients we obtain the solution given in
the corollary.
\end{proof}
\begin{remark}
Using the formula for the Gauss curvature of the surface with $F=0$
\begin{equation}
K=\frac{-1}{\sqrt{EG}}\left(\left(\frac{(\sqrt{G})_u}{\sqrt{E}}\right)_u+\left(\frac{(\sqrt{E})_v}{\sqrt{G}}\right)_v\right).
\end{equation}
we can check that the target surface in Corollary \ref{P1} has
positive Gaussian curvature $K=\frac{a^2}{(a^2-u^2)^2}$, the target
surface in Corollary \ref{Helicoid} has negative Gaussian curvature
$K=\frac{-a^2}{(a^2+u^2)^2}$, and the target surface in Corollary
\ref{Cone} has zero Gaussian curvature $K=0$. Therefore, our results
show that there exist many proper biharmonic maps from a noncompact
surface into surfaces with $K>0$, $K=0$ and $K<0$. This is a
contrast to a theorem proved in \cite{Ji1} asserting that any
biharmonic map from a compact orientable manifold into a manifold
$N$ with sectional curvature ${\rm Rie}^N\leq 0$ has to be harmonic.
\end{remark}
\begin{corollary}
Let $\mathbb{R}^2$ be the Euclidean plane coordinates $(x,y)$, i.e.,
$ds^2=dx^2+dy^2$ and $N =\mathbb{R}\times \mathbb{R}$ with
coordinates $(x,y)$ be provided with the warped product metric
$d\rho^2=dx^2+\sigma^2(x)dy^2$. Then, the identity map ${\bf i}:
\mathbb{R}^2\longrightarrow N$ defined by ${\bf i} \,(x,y)=(x,y)$
 is biharmonic if and only if
\begin{equation}\label{4E}
\begin{cases}
 \sigma^2(x)=-4\ln|x+c_1|;\\
\sigma^2(x)=-4\ln |\cos (\frac{a}{2}x+c_1)|+c_2, {\rm or}\\
\sigma^2(x)=-4b\ln (e^{bx/2}-c_1e^{-bx/2})+c_2.
\end{cases}
\end{equation}
\end{corollary}
\begin{proof}
Applying Proposition \ref{P0} with $f(x)=x$ we conclude that the
identity map ${\bf i}: \mathbb{R}^2\longrightarrow N$ is biharmonic
if and only if
\begin{equation}
(\sigma \sigma')''-(\sigma \sigma')(\sigma \sigma')'=0.
\end{equation}
Setting $y=\sigma \sigma'$ we have
\begin{equation}\label{GD01}
y''-yy'=0,
\end{equation}
which can be written as
\begin{equation}
 (y'-y^2/2)'=0.
\end{equation}
It follows that $\frac{d y}{d x}=(y^2+C)/2$, or
\begin{equation}\label{CMAP}
\frac{d y}{y^2+C}=\frac{1}{2}dx.
\end{equation}
Solving Equation (\ref{CMAP}) we have
\begin{equation}\label{3E}
\begin{cases}
 {\rm For}\;C=0,\;,\;\;\;\;\;\;\;\;\;\;\;\;y=-\frac{2}{x+c_1};\\
{\rm For}\;C=a^2>0.\;\;\;\;\;\;\;y=a\tan (\frac{a}{2}x+c_1);\\
{\rm For}\;C=-b^2<0.\;\;\;\;\;y=\frac
{-b(e^{bx/2}+c_1e^{-bx/2})}{e^{bx/2}-c_1e^{-bx/2}}.
\end{cases}
\end{equation}
Substituting $y=\sigma \sigma'=(\frac{\sigma^2}{2})'$ into
(\ref{3E}) and integrating the resulting equations we obtain
solutions (\ref{4E}). All solutions give proper biharmonic maps as
$\sigma^2\ne {\rm constant}$ in any case. Therefore, we obtain the
corollary.
\end{proof}
\begin{remark}
For more general results on biharmonic identity maps from a product
space into a warped product space see \cite{BMO1}.
\end{remark}

We will end this section with a classification of linear biharmonic
maps between hyperbolic planes with the model
$H^2=(\mathbb{R}^2,e^{-2y}dx^2+dy^2)$.
\begin{proposition}\label{PLast}
Let $\phi:H^2=(\mathbb{R}^2,e^{-2y}dx^2+dy^2)\longrightarrow
H^2=(\mathbb{R}^2,e^{-2v}du^2+dv^2)$ with $\phi(x,y)=(ax+by, cx+dy)$
be a linear map. Then, $\phi$ is biharmonic if and only if $a=b=0$
or $a=\pm1,b=c=0,d=1$; Furthermore, for $a=b=0$ and $d\neq 0$, the
map $\phi$ is a proper biharmonic map.
\end{proposition}
\begin{proof}
One can easily compute the connection coefficients of the domain and
the target surfaces to get
\begin{eqnarray}
&& \Gamma^1_{11}=0, \hskip0.3cm\Gamma^1_{12}=-1, \hskip0.3cm
\Gamma^1_{22}=0,\\\notag && \Gamma^2_{11}=e^{-2y}, \hskip0.3cm
\Gamma^2_{12}=0, \hskip0.3cm \Gamma^2_{22}=0,\\\notag
\end{eqnarray}
and
\begin{eqnarray}
&&\bar{\Gamma}^1_{11}=0, \hskip0.3cm\bar{\Gamma}^1_{12}=-1,
\hskip0.3cm \bar{\Gamma}^1_{22}=0,\\\notag &&
\bar{\Gamma}^2_{11}=e^{-2v}, \hskip0.3cm \bar{\Gamma}^2_{12}=0,
\hskip0.3cm \bar{\Gamma}^2_{22}=0.\\\notag
\end{eqnarray}
We can also check that the components of the Riemannian curvature of
the target surface are given by

\begin{eqnarray}
&&\bar{R}^1_{212}=-1,\; \bar{R}^1_{221}=1,\;
\bar{R}^2_{112}=e^{-2v}, \;\bar{R}^2_{121}=-e^{-2v}\\\notag && {\rm
others}\; \quad \bar{R}^l_{kij}=0,
\end{eqnarray}

and the tension field of the map $\phi$ has components
\begin{eqnarray}\label{TF1}
\tau^1(\phi) &=&-b-2bd-2ace^{2y},\\\label{TF2} \tau^2(\phi)
&=&-d+a^2e^{2y-2v}+b^2e^{-2v}.
\end{eqnarray}

A further computation gives

\begin{eqnarray}\notag
&&\tau^1_1=0,\;\;\tau^1_{11}=0
,\;\;\tau^1_2=-4ace^{2y},\;\;\tau^1_{22}=-8ace^{2y},\\\notag
&&\tau^2_1=-2a^2ce^{2y-2v}-2b^2ce^{-2v},\\\notag
&&\tau^2_{11}=4a^2c^2e^{2y-2v}+4b^2c^2e^{-2v},\\
&&\tau^2_2=2a^2(1-d)e^{2y-2v}-2b^2de^{-2v},\\\notag
&&\tau^2_{22}=4a^2(1-d)^2e^{2y-2v}+4b^2d^2e^{-2v},\\\notag
\end{eqnarray}

\begin{eqnarray}\label{12}
\Delta\tau^1=g^{ij}(\tau^1_{ij}-\Gamma^k_{ij}\tau^1_k)=-4ace^{2y},
\end{eqnarray}

\begin{eqnarray}
&&2g(\nabla\tau^{\alpha},\nabla\phi^{\beta})\bar{\Gamma}^1_{\alpha\beta}\\\notag
&&=4a^3ce^{4y-2v}+4(ab^2c-a^2b(1-d))e^{2y-2v}+4b^3de^{-2v}+8acde^{2y},\\\notag
\end{eqnarray}

\begin{eqnarray}
\tau^{\alpha}\Delta\phi^{\beta}\bar{\Gamma}^1_{\alpha\beta}=a^2be^{2y-2v}+b^3e^{-2v}-2acde^{2y}-2bd-2bd^2,
\end{eqnarray}

\begin{eqnarray}
&&\tau^{\alpha}g(\nabla\phi^{\beta},\nabla\phi^{\rho})(\partial_{\rho}\bar{\Gamma}^1_{\alpha\beta}
+\bar{\Gamma}^v_{\alpha\beta}\bar{\Gamma}^1_{v\rho})\\\notag
=&&\tau^1(g^{11}c^2+d^2)+\tau^2(g^{11}ac+bd)-e^{-2v}\tau^1(g^{11}a^2+b^2)\\\notag
=&&3a^3ce^{4y-2v}+(a^2b+3ab^2c+3a^2bd)e^{2y-2v}-2ac^3e^{4y}\\\notag
&&+(-bc^2-acd-2acd^2-2bc^2d)e^{2y}\\\notag
&&+(b^3+3b^3d)e^{-2v}-2bd^2-2bd^3,
\end{eqnarray}

\begin{eqnarray}
&&-\tau^vg(\nabla\phi^{\alpha},\nabla\phi^{\beta})\bar{R}^1_{\beta\alpha
v}\\\notag
=&&-\tau^1g(\nabla\phi^2,\nabla\phi^2)\bar{R}^1_{221}-\tau^2g(\nabla\phi^1,\nabla\phi^2)\bar{R}^1_{212}\\\notag
=&&2ac^3e^{4y}+(bc^2+2bc^2d+2acd^2-acd)e^{2y}+a^3ce^{4y-2v}\\\notag
&&+(a^2bd+ab^2c)e^{2y-2v}+b^3de^{-2v}+2bd^3.\\\notag
\end{eqnarray}

Similarly, we have
\begin{eqnarray}
\Delta\tau^2&=&g^{ij}(\tau^2_{ij}-\Gamma^k_{ij}\tau^2_k)\\\notag
&=&(2b^2d+4b^2d^2)e^{-2v}+(2a^2+4b^2c^2-6a^2d+4a^2d^2)e^{2y-2v}\\\notag
&&+4a^2c^2e^{4y-2v},
\end{eqnarray}
\begin{eqnarray}
2g(\nabla\tau^{\alpha},\nabla\phi^{\beta})\bar{\Gamma}^2_{\alpha\beta}=2g(\nabla\tau^1,\nabla\phi^1)\bar{\Gamma}^2_{11}=-8abce^{2y-2v},
\end{eqnarray}

\begin{eqnarray}
\tau^{\alpha}\Delta\phi^{\beta}\bar{\Gamma}^2_{\alpha\beta}=-b\tau^1\bar{\Gamma}^2_{11}=(b^2+2b^2d)e^{-2v}+2abce^{2y-2v},
\end{eqnarray}

\begin{eqnarray}
&&\tau^{\alpha}g(\nabla\phi^{\beta},\nabla\phi^{\rho})(\partial_{\rho}\bar{\Gamma}^2_{\alpha\beta}
+\bar{\Gamma}^v_{\alpha\beta}\bar{\Gamma}^2_{v\rho})\\\notag
=&&-2\tau^1(g^{11}ac+g^{22}bd)e^{-2v}+\tau^1(g^{11}ac+g^{22}bd)\bar\Gamma^1_{12}\bar\Gamma^2_{11}\\\notag
&&+\tau^2(g^{11}a^2+g^{22}b^2)\bar\Gamma^1_{21}\bar\Gamma^2_{11}\\\notag
=&&-b^4e^{-4v}+(4b^2d+6b^2d^2)e^{-2v}-2a^2b^2e^{2y-4v}\\\notag
&&+(3abc+a^2d+12abcd)e^{2y-2v}-a^4e^{4y-4v}+6a^2c^2e^{4y-2v},\\\notag
\end{eqnarray}
and
\begin{eqnarray}\label{13}
&&-\tau^vg(\nabla\phi^{\alpha},\nabla\phi^{\beta})\bar{R}^2_{\beta\alpha
v}\\\notag
=&&-\tau^1g(\nabla\phi^2,\nabla\phi^1)\bar{R}^2_{121}-\tau^2g(\nabla\phi^1,\nabla\phi^1)\bar{R}^2_{112}\\\notag
=&&(a^2d-abc-4abcd)e^{2y-2v}-a^4e^{4y-4v}-2a^2b^2e^{2y-4v}\\\notag
&&-b^4e^{-4v}-2b^2d^2e^{-2v}-2a^2c^2e^{4y-2v}.\\\notag
\end{eqnarray}
Substitute (\ref{12})$\sim$ (\ref{13}) into (\ref{BI3}), we conclude
that $\phi$ is biharmonic if and only if it solves the system
\begin{equation}\label{14}
\begin{cases}
(-4ac+4acd)e^{2y}+8a^3ce^{4y-2v}+(8ab^2c-2a^2b+8a^2bd)e^{2y-2v}\\
+(8b^3d+2b^3)e^{-2v}-2bd-4bd^2=0,\\
\\
(8b^2d+8b^2d^2+b^2)e^{-2v}+8a^2c^2e^{4y-2v}-2b^4e^{-4v}-4a^2b^2e^{2y-4v}\\
-2a^4e^{4y-4v}+(2a^2+4b^2c^2-4a^2d+4a^2d^2-4abc+8abcd)e^{2y-2v}=0.
\end{cases}
\end{equation}
The first equation of (\ref{14}) is equivalent to
\begin{eqnarray}\label{15}
&&e^{-2cx}[8a^3ce^{4y-2dy}+(8ab^2c-2a^2b+8a^2bd)e^{2y-2dy}\\\notag
&&+(8b^3d+2b^3)e^{-2dy}]+(-4ac+4acd)e^{2y}-2bd-4bd^2=0
\end{eqnarray}
for any $x,y\in R$.\\
If $c\neq 0$, (\ref{15}) implies\\
\begin{equation*}\label{GE}
\begin{cases}
8a^3ce^{4y-2dy}+(8ab^2c-2a^2b+8a^2bd)e^{2y-2dy}+(8b^3d+2b^3)e^{-2dy}=0,\\
\\
-4ac+4acd=0,\\
\\
2bd-4bd^2=0.
\end{cases}
\end{equation*}
Solving this equation we obtain $a=b=0,$ $c\neq 0$.\\
\\
If $c=0$, (\ref{15}) becomes \\
\begin{eqnarray}\label{16}
(-2a^2b+8a^2bd)e^{2(1-d)y}+(8b^3d+2b^3)e^{-2dy}-2bd-4bd^2=0,
\end{eqnarray}
which can be solved by considering the following three cases:\\
(A) For $d=0$, (\ref{16}) becomes $-2a^2be^{2y}+2b^3=0$, which has solution $b=d=c=0$;\\
(B) For $d=1$, (\ref{16}) becomes $10b^3e^{-2y}+6a^2b-6b=0$ whose
solution is  $b=c=0,$ $d=1$;\\
(C) For $d\neq 0$, and $d\neq 1$, (\ref{16}) implies that\\
\begin{equation*}\label{GE}
\begin{cases}
-2a^2b+8a^2bd=0,\\
\\
8b^3d+2b^3=0,\\
\\
-2bd-4bd^2=0.
\end{cases}
\end{equation*}
Solving this we have $b=c=0$, $d\neq0,\ d\neq1$.\\

So, the solution of the 1st equation of (\ref{14}) is $(i)$ $a=b=0$;
and
$(ii)$ $b=c=0$.\\
\\
It is easy to check that $a=b=0$ is also a solution of the 2nd
equation of (\ref{14}), and hence $a=b=0$ is a solution of (\ref{14}).\\

Substituting $b=c=0$ into the 2nd equation of (\ref{14}) we obtain
$(2a^2-4a^2d+4a^2d^2)e^{2(1-d)y}-2a^4e^{4(1-d)y}=0$. Solving this,
we have $a=b=c=0,\; d\neq 1$, or $a=b=c=0, \;d=1$, or $a=\pm1, \;b=c=0, \;d=1$.\\

Summarizing the above results we conclude that the Equation
(\ref{14}) has solutions $a=b=0$, or $a=\pm1, \;b=c=0, \;d=1$.\\

Using the components (\ref{TF1}) and (\ref{TF2}) of the tension
field of $\phi$ we can check that when $a=b=0,\;d\neq0$, the map
$\phi(x,y)=(ax+by,cx+dy)$ is a proper biharmonic map. This completes
the proof of the proposition.
\end{proof}

\begin{remark}
Note that Oniciuc \cite{On} proves that if $\phi : (M, g)
\longrightarrow (N, h)$ is a map with the property that
$|\tau(\phi)| =$ constant,  ${\rm Riem}^N < 0$, and there is a point
$p \in M$ such that $\rank_p\, \phi\ge 2$, then $\phi$ is biharmonic
if and only if it is harmonic. Our Proposition \ref{PLast} shows
that there does exist a lot of proper biharmonic maps into a space
with negative curvature, having the property $|\tau(\phi)| = {\rm
constant}$ but $\rank_p\, \phi< 2$ for every point $p$.
\end{remark}

\end{document}